\documentclass[12pt,a4]{article}
\usepackage{indentfirst}
\usepackage{amssymb}
\newtheorem{theo}{Theorem}
\newtheorem{p}[theo]{Proposition}

\newtheorem{de}[theo]{Definition}

\newcommand{\pr}{Proposition }

\newcommand{\qed}{\hfill {\rule{3mm}{3mm}}}
\newcommand{\ma}[5]{\ensuremath{#1:#2\longrightarrow #3\,, \quad #4 \longmapsto  #5}}
\newcommand{\map}[4]{\ensuremath{#1\longrightarrow #2, \quad #3 \longmapsto   #4}}
\newcommand{\mb}[3]{\ensuremath{#1:#2\longrightarrow #3}}
\newcommand{\me}[2]{\ensuremath{\left\{\,\left. #1\;\right|\; #2 \, \right\}}}

\newcommand{\sa}[2]{\ensuremath{\left \langle \, #1 \,, \, #2 \, \right \rangle}}

\newcommand{\bk}{\ensuremath{\mathrm{I\! K}}}
\newcommand{\n}[1]{\ensuremath{\left\| #1 \right\|}}
\renewcommand{\labelenumi}{\alph{enumi})}
\begin{document}

\begin{center}
\Large {\bfseries
Uniqueness of W*-tensor Products} \vspace{0.5cm}\\ 

\end{center}
\begin{flushright}
\large \bfseries
Corneliu Constantinescu
\end{flushright}
\vspace{1cm}
\begin{abstract}
In contrast to C*-algebras, distinct C*-norms on the algebraic tensor product of two W*-algebras produce isomorphic W*-tensor products.  
\end{abstract}

AMS classification code: 46L06, 46L10
 

We use the notation and terminology of [C]. For W*-tensor products of W*-algebras we use [T].

In the sequel we give a list of some notation used in this paper.
\renewcommand{\labelenumi}{\arabic{enumi}.}
\begin{enumerate}
\item  $\bk$ denotes the field of real or the field of complex numbers. The whole theory is developed in parallel for the real and complex case (but the proofs coincide).
\item If $f$ is a map defined on a set $X$ and $Y$ is a subset of $X$ then $f|Y$ denotes the restriction of $f$ to $Y$.
\item If $E,F$ are vector spaces in duality then $E_F$ denotes the vector space $E$ endowed with the locally convex topology of pointwise convergence on $F$, i.e. with the weak topology $\sigma (E,F)$.
\item If $E$ is a Banach space then $E'$ denotes its dual, $E''$ its bi-dual, and  $E^{\#}$ its unit ball :
$$ E^{\#}:= \me{x\in E}{\|  x\|  \leq 1}.$$
We put for every $x\in E$, 
$$\ma{j_Ex}{E'}{\bk}{x'}{\sa{x}{x'}}$$
and call the map $\mb{j_E}{E}{E''}$ the evaluation map of $E$.
If $F$ is a vector subspace of $E$ then we set
$$F^0:=\me{x'\in E'}{x'|F=0}. $$
If $F$ is a vector subspace of $E'$ then we define
$$^0F:=\me{x\in E}{x'\in F\Longrightarrow \sa{x}{x'}=0}.$$
\item Let $E,F$ be Banach spaces and $\mb{\varphi }{E}{F}$ is a continuous linear map. $\varphi $ is called an isometry if it preserves the norms and if it is surjective. We put $Im\,\varphi :=\varphi (E)=\me{\varphi x}{x\in E}$  and denote by 
$$\ma{\varphi '}{F'}{E'}{y'}{y'\circ \varphi }$$
the transpose of $\varphi $. 
\item If $E$ is a C*-algebra then we denote by $Pr\,E$ the set of orthogonal projections of $E$. If in addition $E$ is unital then we denote by $1_E$ its unit. If $F$ is a subset of $E$ then we put
$$F^c:=\me{x\in E}{y\in F\Longrightarrow xy=yx}.$$
We put for all $(x,x')\in E\times E'$ 
$$\ma{xx'}{E}{\bk}{y}{\sa{yx}{x'}},$$
$$\ma{x'x}{E}{\bk}{y}{\sa{xy}{x'}};$$
then $xx',x'x\in E'$.
\item If $E$ is a W*-algebra then $\ddot{E} $ denotes its predual.
\item $\odot$ denotes the algebraic tensor product. If $E,F$ are C*-algebras and if $\alpha $ is a C*-norm on $E\odot F$ then $E\otimes _\alpha F$ denotes the C*-algebra obtained by the compltion of $E\odot F$ with respect to $\alpha $.
\end{enumerate}

\renewcommand{\labelenumi}{\alph{enumi})} 

\begin{p}\label{8085}
Let $E$ be a W*-algebra and $F$ a closed vector subspace of $\ddot{E} $ such that $xF\subset F$ and $Fx\subset F$ for all $x\in E$.
\begin{enumerate}
\item There is a $p\in Pr\, E^c$ such that $F=p\ddot{E} $ and $F^0=(1_E-p)E$.
\item For every $x\in pE$ put
$$\ma{\tilde{x} }{F}{\bk}{a}{\sa{x}{a}}.$$
Then $\tilde{x}\in F' $ for all $x\in pE$ and the map
$$\map{pE}{F'}{x}{\tilde{x} }$$
is an isometry of Banach spaces.
\end{enumerate}
\end{p}

a) follows from [T] Theorem III.2.7 c).

b) For $a\in \ddot{E} $,
$$\sa{x}{a}=\sa{px}{a}=\sa{x}{ap}=\sa{\tilde{x} }{ap},$$
so $\tilde{x}\in F' $ and $\n{x}=\n{\tilde{x} }$. Let $a'\in F'$ and put
$$\ma{y}{\ddot{E} }{\bk}{a}{\sa{a'}{pa}}.$$
Then $y\in E$ and for $a\in \ddot{E} $,
$$\sa{py}{a}=\sa{y}{ap}=\sa{y}{a},$$
so $y=py\in pE$ and $\tilde{y}=a' $, i.e. the map is surjective.\qed

\begin{p}\label{8086}
Let $E$ be a C*-algebra and $F$ a closed vector subspace of $E'$ such that $xF\subset F$ and $Fx\subset F$ for all $x\in E$.
\begin{enumerate}
\item There is a $p\in Pr\,(E'')^c$ such that $F=pE'$ and $F^0=(1_{E''}-p)E''$.
\item The map 
$$\map{pE''}{F'}{x''}{x''|F}$$
is an isometry of Banach spaces.
\end{enumerate}
\end{p}

By [C] Corollary 1.3.6.5, $Im\,j_E$ is dense in $E''_{E'}$ so $x''F\subset F$ and $Fx''\subset F$ for all $x''\in E''$. By [C] Theorem 6.3.2.1 b), $E''$ is a W*-algebra and the assertions follow from \pr\ref{8085}.\qed
 
\begin{de}\label{8088}
Let $E,F$ be W*-algebras. We define the (bilinear) maps
$$\map{(E\odot F)\times (\ddot{E}\odot \ddot{F}  )}{\ddot{E}\odot \ddot{F}  }{(x\otimes y,a\otimes b)}{(x\otimes y)(a\otimes b):=(xa)\otimes (yb)},$$
$$\map{(\ddot{E}\odot \ddot{F}  )\times (E\odot F)}{\ddot{E}\odot \ddot{F}  }{(a\otimes b,x\otimes y)}{(a\otimes b)(x\otimes y):=(ax)\otimes (by)}.$$
and similarly the (bilinear) maps
$$(E\odot F)\times (E'\odot F')\longrightarrow E'\odot F',$$
$$(E'\odot F')\times (E\odot F)\longrightarrow E'\odot F'.$$
\end{de}

\begin{p}\label{8087}
Let $E,F$ be W*-algebras, $\alpha $ a C*-norm on $E\odot F$, and 
$$G:=\ddot{E}\bar{\otimes }_\alpha \ddot{F}$$
the closure of $Im\,j_{\ddot{E} }\odot Im\,j_{\ddot{F} }$ in $(E\otimes _\alpha F)'$.
\begin{enumerate}
\item There is a $p\in Pr\,((E\otimes _\alpha F)'')^c$ such that $G=p(E\otimes _\alpha F)'$ and such that the map
$$\map{p(E\otimes _\alpha F)''}{G'}{x''}{x''|G}$$
is an isometry of Banach spaces. 
\item If $\mb{j}{E\otimes _\alpha F}{(E\otimes _\alpha F)''}$ denotes the evaluation map then $Im\,j$ is a C*-subalgebra of $p(E\otimes _\alpha F)''$ generating it as a W*-algebra.
\end{enumerate}
\end{p}

a) For $(u,v),(x,y)\in E\times F$ and $(a,b)\in \ddot{E}\times \ddot{F}  $,
$$\sa{u\otimes v}{(x\otimes y)((j_{\ddot{E} }a)\otimes (j_{\ddot{F}}b))}=\sa{(u\otimes v)(x\otimes y)}{(j_{\ddot{E}}a)\otimes (j_{\ddot{F} }b)}=$$
$$=\sa{(ux)\otimes (vy)}{(j_{\ddot{E}}a)\otimes (j_{\ddot{F} }b)}=\sa{ux}{j_{\ddot{E} }a}\sa{vy}{j_{\ddot{F} }b}=$$
$$=\sa{ux}{a}\sa{vy}{b}=\sa{u}{xa}\sa{v}{yb}
=\sa{u}{j_{\ddot{E} }(xa)}\sa{v}{j_{\ddot{F} }(yb)}=$$
$$=\sa{u\otimes v}{j_{\ddot{E}}(xa)\otimes j_{\ddot{F} }(yb)},$$
so
$$(x\otimes y)((j_{\ddot{E}}a)\otimes (j_{\ddot{F} }b))=j_{\ddot{E}}(xa)\otimes j_{\ddot{F} }(yb)\in G.$$
It follows $(x\otimes y)G\subset G$ and $zG\subset G$ for all $z\in E\otimes _\alpha F$. Similarly $Gz\subset G$ for all $z\in E\otimes _\alpha F$. By  \pr\ref{8086}, there is a $p\in Pr\,((E\otimes _\alpha F)'')^c$ such that $G=p(E\otimes _\alpha F)'$ and such that the map
$$\map{p(E\otimes _\alpha F)''}{G'}{x''}{x''|G}$$
is an isometry of Banach spaces.

b Let $(x,y)\in E\times F$ and $\varepsilon >0$. There are $a\in E^{\#}$ and $b\in F^{\#}$ such that
$$\n{x}\n{y}-\varepsilon <\sa{x}{a}\sa{y}{b}=\sa{x}{j_{\ddot{E} }a}\sa{y}{j_{\ddot{F} }b}=$$
$$=\sa{x\otimes y}{(j_{\ddot{E} }a)\otimes (j_{\ddot{F} }b)}=\sa{j(x\otimes y)}{(j_{\ddot{E} }a)\otimes (j_{\ddot{F} }b)}\leq \n{j(x\otimes y)|G}.$$
Since $\varepsilon $ is arbitrary,
$$\n{j(x\otimes y)|G}=\n{x}\n{y},$$
so by a), $j(x\otimes y)\in p(E\otimes _\alpha F)''$. It follows $j(E\odot F)\subset p(E\otimes _\alpha F)''$ and $j(E\otimes _\alpha  F)\subset p(E\otimes _\alpha F)''$. Let 
$$z\in (Im\,j_{\ddot{E}}\odot Im\,j_{\ddot{F} })\cap \,^0(j(E\otimes _\alpha F)).$$
Then $z|(E\otimes _\alpha F)=0$ and so $z=0$. Thus
$$G\cap \,^0(j(E\otimes _\alpha F))=\{0\},$$
so by a),
$$^0(j(E\otimes _\alpha F))\subset (1_{(E\otimes _\alpha F)''}-p)(E\otimes _\alpha F)',$$
$$p(E\otimes _\alpha F)''\subset (\,^0(j(E\otimes _\alpha F)))^0,$$
$$p(E\otimes _\alpha F)''= (\,^0(j(E\otimes _\alpha F)))^0.$$
Hence $p(E\otimes _\alpha F)''$ is the closure of $j(E\otimes _\alpha F)$ in $(E\otimes _\alpha F)''_{(E\otimes_ \alpha F)'}$ ([C] Proposition 1.3.5.4). By [C] Corollary 4.4.4.12 a), $p(E\otimes _\alpha F)''$ is the W*-subalgebra of $(E\otimes _\alpha F)''$ generated by $j(E\otimes _\alpha F)$.\qed

\begin{de}\label{8090}
We put \emph{(with the notation of \pr\ref{8087})} 
$$E\bar{\otimes }_\alpha F:=p(E\otimes _\alpha F)'' .$$ 
It is a W*-algebra with $\ddot{E}\bar{\otimes }_\alpha \ddot{F}$ as predual.
\end{de}

\begin{theo}\label{8091}
If $E,F$ are W*-algebras and $\alpha ,\beta $ are C*-norms on $E\odot F$ then $E\bar{\otimes }_\alpha F$ and $E\bar{\otimes }_\beta  F$ are isomorphic.
\end{theo}

We may assume $\alpha \leq \beta $. By [W] Proposition T.6.24, there is a surjective C*-homomorphism 
$$\mb{\varphi }{E\otimes _\beta F}{E\otimes _\alpha F}.$$ 
Then
$$\mb{\varphi '}{(E\otimes _\alpha F)'}{(E\otimes _\beta F)'}$$
preserves the norms ([C] Proposition 1.3.5.2). Thus $\ddot{E}\bar{\otimes }_\alpha \ddot{F}=\ddot{E}\bar{\otimes }_\beta  \ddot{F}$ and $\varphi ''$ is an isometry of Banach spaces. By [C] Corollary 6.3.2.3, $\varphi ''$ is a W*-isomorphism.\qed

\begin{center}
{\bfseries REFERENCES}
\end{center}
\begin{flushleft}
[C] Constantinescu, Corneliu, C*-algebras. Elsevir, 2001. \newline
[T] Takesaki, Masamichi, Theory of Operator Algebra I. Springer, 2002. \newline
[W] Wegge-Olsen, N. E., K-theory and C*-algebras. Oxford University Press, 1993. \newline
\end{flushleft}
\begin{flushright}
{\scriptsize \hspace{-5mm} Corneliu Constantinescu$\quad$\\
Bodenacherstr. 53$\qquad\;$\\
CH 8121 Benglen$\qquad\;\;$\\
e-mail: constant@math.ethz.ch }
\end{flushright}
\end{document}